\documentclass[12pt]{amsart}
\usepackage{amssymb, euscript}
\input diagrams.sty
\makeatletter
\@addtoreset{equation}{section}
\makeatother

\newcommand{\xref}[1]{{\rm \ref{#1}}}

\newcommand{\mt}[1]{\operatorname{#1}}

\newcommand{\Bs}{\operatorname{Bs}}
\newcommand{\mult}{\operatorname{mult}}

\renewcommand{\emptyset}{\varnothing}
\newcommand{\var}{\varphi}

\newcommand{\CC}{{\mathbb C}}

\newcommand{\ZZ}{{\mathbb Z}}
\newcommand{\QQ}{{\mathbb Q}}
\newcommand{\PP}{{\mathbb P}}

\newcommand{\NN}{{\mathbb N}}
\newcommand{\WW}{{\mathbb W}}
\newcommand{\OOO}{{\mathcal O}}

\newcommand{\LLL}{{\EuScript{L}}}

\newcommand{\HHH}{{\EuScript{H}}}
\newcommand{\EEE}{{\EuScript{E}}}

\newcommand{\ov}[1]{\overline{#1}}

\newtheorem{theorem}[equation]{Theorem}
\newtheorem{claim}[equation]{Claim}

\newtheorem{proposition}[equation]{Proposition}
\newtheorem{proposition-definition}[equation]{Proposition-Definition}
\newtheorem{lemma}[equation]{Lemma}
\newtheorem{corollary}[equation]{Corollary}

\theoremstyle{definition}

\newtheorem{example}[equation]{Example}
\newtheorem{examples}[equation]{Examples}
\newtheorem{emptytheorem}[equation]{}
\newtheorem{remark}[equation]{Remark}

\newcounter{ssylki}[equation]
\renewcommand{\thessylki}{(\arabic{section}.%
\arabic{equation}.\arabic{ssylki})}
\newenvironment{subth}{\refstepcounter{ssylki}
\par\medskip\noindent{\rm\thessylki}}{\par}

\title{A remark on Fano threefolds\\
with canonical Gorenstein singularities}
\author{Yuri G. Prokhorov}
\thanks
{This work was is carried out under the support of the grant
INTAS-OPEN 2000-269}

\address{
Department of Algebra, Faculty of Mathematics, Moscow State
Lomonosov University, Leninskie Gory, 117234 Moscow, Russia}
\address{
Research Institute for Mathematical Sciences, Kyoto University,
Kyoto 606-8502, Japan}

\email{prokhoro@mech.math.msu.su}

\begin{document}
\begin{abstract}
We give some rationality constructions for Fano threefolds with
canonical Gorenstein singularities.
\end{abstract}

\maketitle
\section{Introduction}
A normal projective variety $X$ over $\CC$ is called \emph{Fano
variety} if its anticanonical Weil divisor $-K_X$ is Cartier and
ample. The number $(-K_X)^{\dim X}$ is called the \emph{degree} of
$X$.

In this paper we deal with Fano threefolds having at worst
canonical Gorenstein singularities. The motivation to consider
such Fanos is the following result due to Alexeev:

\begin{theorem}[\cite{A}]
\label{th-alexeev}
Let $Y$ be a $\QQ$-Fano threefold with terminal $\QQ$-factorial
singularities and Picard number $1$. If the anticanonical model
$\Phi_{|-K_Y|}(Y)$ is three-dimensional, then $Y$ is birationally
equivalent to a Fano threefold $X$ with Gorenstein canonical
singularities and base point free $|-K_{X}|$.
\end{theorem}

Smooth Fano threefolds were classified by Iskovskikh and
Mori-Mukai. If a Fano threefold $X$ has only terminal Gorenstein
singularities, then it has a smoothing (see \cite{Nam}), i.e.,
such $X$ can be considered as a degeneration of a smooth ones. One
can expect the same situation in the case of cDV singularities. In
contrast, Fano threefolds with canonical non-cDV singularities are
not necessarily smoothable:

\begin{examples}
\label{ex-Fano}
(i) (Weighted projective spaces.) All weighted projective spaces
$\PP(a_1,\dots,a_n)$ are $\QQ$-factorial $\QQ$-Fanos with only log
terminal singularities and $\rho=1$. It is easy to enumerate all
these $3$-dimensional spaces (up to isomorphisms) having at worst
Gorenstein canonical singularities: $\PP(4,3,3,2)$,
$\PP(15,10,3,2)$, $\PP(4,4,3,1)$, $\PP(10,5,4,1)$,
$\PP(21,14,6,1)$, $\PP(6,3,2,1)$, $\PP(12,8,3,1)$, $\PP(5,2,2,1)$,
$\PP(2,2,1,1)$, $\PP(9,6,2,1)$, $\PP(4,2,1,1)$,
$\PP(1,1,1,1)=\PP^3$, $\PP(6,4,1,1)$, $\PP(3,1,1,1)$. For
$X=\PP(6,4,1,1)$ and $\PP(3,1,1,1)$ we have $-K_X^3=72$ while
there are no smooth Fano threefolds of degree $72$.

(ii) (Cones) Let $S$ be a del Pezzo surface of degree $d$ with at
worst Du Val singularities and let $\LLL=\OOO_S(-K_S)$. Consider
the $\PP^1$-bundle $\PP=\PP_S(\OOO_S\oplus \LLL)$. The map $\PP\to
X$ given by the linear system $|\OOO_{\PP}(n)|$, $n\gg 0$
contracts the negative section. Since $-K_{\PP}\sim
\OOO_{\PP}(2)$, the variety $X$ is a Fano threefold of index $2$
and degree $8d$ with canonical Gorenstein singularities. For
$S=\PP^2$ we have $-K_X^3=72$ and $X\simeq \PP(3,1,1,1)$.
\end{examples}

In this paper we propose a rationality construction for Fano
threefolds such as in Theorem \xref{th-alexeev}. Our main result
is the following.

\begin{theorem}
\label{rational}
Let $X$ be a Fano threefold of degree $-K_X^3=2g-2$ with at worst
canonical Gorenstein singularities. Assume that the linear system
$|-K_X|$ is base point free and $X$ has at least one non-cDV
point. If $X$ is not hyperelliptic, then $X$ is rational except
for two cases:

\begin{subth}
\label{case-quartic}
$X$ is a quartic in $\PP^4$, or
\end{subth}

\begin{subth}
\label{case-quadric-cubic}
$X$ is an intersection of a quadric and a cubic in $\PP^5$.
\end{subth}

\par\medskip\noindent
If $X$ is hyperelliptic, then $X$ is rational except perhaps for
the following cases:
\begin{subth}
\label{case-hyperell-conic}
the anticanonical image of $X$ is a projective cone over a surface
degree $g-1$ in $\PP^g$, and $X$ is birationally equivalent to a
conic bundle, or
\end{subth}
\begin{subth}
\label{case-hyperell-other}
$g\le 4$ and $X$ is birationally equivalent to one of varieties
\xref{gather-cases-1}-\xref{gather-cases-4} below.
\end{subth}
\end{theorem}

\begin{remark}
(i) There are examples of varieties as in \xref{case-quartic} and
\xref{case-quadric-cubic} which are rational as well as
nonrational (see \xref{ex-quartic} and \xref{ex-quardic-cubic}).
The answer to the rationality question in these cases depends on
delicate analysis of singular points.

(ii) Biregular theory of Fano threefolds with canonical Gorenstein
singularities was developed by Mukai \cite{Mu}.
\end{remark}

For the proof we study the projection from a non-cDV point using
technique of \cite{CF} and \cite{A} (cf. \cite{Ch1}).

\subsection*{Acknowledgments}
The work has been completed during my stay at Kyoto Research
Institute for Mathematical Sciences. I would like to thank RIMS
hospitality and support.

\section{Preliminaries}
\subsection*{Notation}
All varieties are assumed to be defined over $\CC$. Let $\HHH$ be
a linear system of Weil divisors on a normal variety.
$\Phi_{\HHH}\colon X\dashrightarrow \PP^{\dim \HHH}$ denotes the
corresponding rational map. By $\PP(a_1,\dots,a_n)$, $a_i\in \NN$
we denote the weighted projective space that is $\mt{Proj}
\CC[x_1,\dots,x_n]$ with grading $\deg x_i=a_i$. Let $\EEE$ be a
locally free sheaf on $X$. Then $\PP_X(\EEE)$ is the
projectivization $\mt{Proj}\EEE^*$. The Picard number of a variety
$V$ is denoted by $\rho(V)$.

\subsection*{$\QQ$-Fano threefolds of index $>1$}
Recall that a Fano threefold $X$ is called a \emph{del Pezzo
threefold} if $X\not\simeq\PP^3$ and $-K_X= 2H$ for some ample
Cartier divisor $H$. The number $H^3$ is called the \emph{degree}
of $X$. Everywhere below we will assume that all del Pezzo
threefolds have at worst canonical Gorenstein singularities.

\begin{theorem}[\cite{CF}]
\label{th-CF}
Let $W$ be a $\QQ$-factorial threefold with at worst terminal
singularities and $\rho(X)=1$ and let $H\subset W$ be a smooth del
Pezzo surface contained into the smooth locus of $X$. Assume that
$W$ is not rational. Then $(W,H)$ is one of the following:
\begin{subth}
$W=W_6\subset \PP(1,1,2,2,3)$ is a weighted hypersurface of degree
$6$, $H\in |\OOO(2)|$,\ $\dim |H|=4$,
\label{gather-cases-1}
\end{subth}
\begin{subth}
$W$ is a smooth cubic in $\PP^4$,\ $H\in |-\frac12 K_W|$,\ $\dim
|H|=4$,
\label{gather-cases-2}
\end{subth}
\begin{subth}
$W$ is a del Pezzo threefold of degree $2$,\ $H\in |-\frac12
K_W|$,\ $\dim |H|=3$,
\label{gather-cases-3}
\end{subth}
\begin{subth}
$W$ is a del Pezzo threefold of degree $1$,\ $H\in |-\frac12
K_W|$,\ $\dim |H|=2$.
\label{gather-cases-4}
\end{subth}
\end{theorem}

It is known that general (and even all smooth) members of families
\xref{gather-cases-3}, \xref{gather-cases-4}, and all members of
the family \xref{gather-cases-2} are nonrational (see \cite{Be}
and \cite{CG}).

\subsection*{Canonical singularities}
A normal three-dimensional singularity $X\ni o$ is said to be
\emph{cDV} if it is analytically isomorphic to a hypersurface
singularity $f(x,y,z)+tg(x,y,z,t)=0$, where $f(x,y,z)=0$ is an
equation of a Du Val singularity. All cDV singularities are
canonical (and Gorenstein).

\begin{theorem}[{\cite{Reid}}]
\label{th-can-sing}
If $X\ni o$ is a cDV point, then the discrepancy of every prime
divisor with center at $o$ is strictly positive. Conversely, if
$X\ni o$ is a canonical three-dimensional singularity such that
the discrepancy of every prime divisor with center at $o$ is
strictly positive, then $X\ni o$ is a cDV point. In particular,
any threefold with canonical singularities has only finitely many
non-cDV points.
\end{theorem}

\subsection*{Singularities of linear systems}
Let $X$ be a normal variety and let $\HHH$ be a movable linear
system (of Weil divisors) on $X$. Assume that $K_X+H$ is
$\QQ$-Cartier for $H\in \HHH$. For any good (but fixed) resolution
$f\colon Y\to X$ of the pair $(X,\HHH)$ we can write
\[
K_Y+\HHH_Y=f^*(K_X+\HHH)+\sum_E a(E,\HHH)E,
\]
where $E$ in the sum runs through all exceptional divisors,
$\HHH_Y$ is the birational transform of $\HHH$, and $a(E,\HHH)\in
\QQ$. Here and below in numerical formulas we easily write linear
system $\HHH$, $\HHH_Y$ etc. instead of their members. We say that
$(X,\HHH)$ has \emph{terminal} (resp. \emph{canonical})
singularities if $a(E,\HHH)>0$ (resp. $a(E,\HHH)\ge 0$) for all
$E$.

\begin{proposition}[\cite{A}]
\label{prop-A-sing-term}
If the pair $(X,\HHH)$ is terminal and members of $\HHH$ are
$\QQ$-Cartier, then $\HHH$ has at worst isolated base points $P_i$
such that $\mult_{P_i}\HHH=1$. In particular, $X$ is smooth at
$P_i$ and $\HHH$ is a linear system of Cartier divisors.
\end{proposition}

In the category of three-dimensional terminal $\QQ$-factorial
pairs the log minimal model program works \cite{A}. In particular,
we have the following.
\begin{proposition-definition}
Let $X$ be a normal threefold and let $\HHH$ be a movable linear
system on $X$. Assume that $K_X+\HHH$ is $\QQ$-Cartier. Then there
is a birational contraction $f\colon Y\to X$ such that
$(Y,\HHH_Y)$ has only terminal singularities and $K_Y+\HHH_Y$ is
$f$-nef, where $\HHH_Y$ is the birational transform of $\HHH$.
Such $f$ is called a \emph{relative terminal model} of $(X,\HHH)$.
In this situation we can write
\[
K_Y+\HHH_Y=f^*(K_X+\HHH)-\sum a_iE_i,
\]
where $a_i\ge 0$ and the $E_i$ are exceptional divisors.
\end{proposition-definition}

\subsection*{Anticanonical linear system}
By Riemann-Roch we have the following.
\begin{lemma}
Let $X$ be a Fano threefold with at worst canonical Gorenstein
singularities. Then $\dim |-K_X|=-\frac12K_X^3+2$.
\end{lemma}
Denote $g=-\frac12K_X^3+1$. This integer is called the
\emph{genus} of $X$. Thus
\[
-K_X^3=2g-2\quad \text{and}\quad \dim |-K_X|=g+1.
\]

\begin{theorem}[\cite{R}]
Let $X$ be a Fano threefold with at worst canonical Gorenstein
singularities. Then the pair $(X,|-K_X|)$ has only canonical
singularities.
\end{theorem}

\subsection*{Anticanonical models}
Similar to the nonsingular case (see \cite{Isk}) one can prove the
following.
\begin{proposition}
\label{prop-hyperell}
Let $X$ be a Fano threefold of genus $g$ with at worst canonical
Gorenstein singularities and let $\Phi=\Phi_{|-K_X|}\colon X
\dashrightarrow \ov X\subset \PP^{g+1}$ be the anticanonical map,
where $\ov X=\Phi(X)$. Then $\dim \ov X\ge 2$ and $\dim \ov X=2$
if and only if $\Bs |-K_X|\neq \emptyset$. If $\Bs
|-K_X|=\emptyset$, then one of the following holds:
\begin{enumerate}
\item
$\Phi\colon X\to \ov X$ is a double covering, in this case $\ov
X\subset \PP^{g+1}$ is a variety of degree $g-1$, or
\item
$\Phi\colon X\to \ov X$ is an isomorphism.
\end{enumerate}
\end{proposition}

\begin{remark}
\label{rem-Enriques}
In case (i) the Fano variety $X$ is called \emph{hyperelliptic}.
Here $\ov X\subset \PP^{g+1}$ is a so-called \emph{variety of
minimal degree}. According to the well-known theorem of Enriques
$\ov X$ is one of the following:
\begin{enumerate}
\item
the image of a $\PP^2$-bundle $\PP_{\PP^1}(\OOO_{\PP^1}(a_1)
\oplus \OOO_{\PP^1}(a_2) \oplus \OOO_{\PP^1}(a_3))$, where $a_i\ge
0$, $\sum a_i>0$ under the morphism defined by the linear system
$|\OOO(1)|$, $g=a_1+a_2+a_3+1$,
\item
$\ov X\subset \PP^4$ is a smooth quadric, $g=3$,
\item
$\ov X\subset \PP^6$ is a cone over the Veronese surface, $g=5$.
\end{enumerate}
\end{remark}

\section{Proof of Theorem \xref{rational}}
Let $X$ be a Fano threefold with only canonical Gorenstein
singularities of genus $g$ such that the anticanonical linear
system $|-K_X|$ is base point free and let
$\Phi=\Phi_{|-K_X|}\colon X \to \ov X\subset \PP^{g+1}$ be the
morphism defined by the anticanonical linear system. Assume that
$X$ has at least one non-cDV point $o\in X$. Let
$\HHH=\HHH_o\subset |-K_X|$ be the subsystem of all divisors
containing $o$. It is clear that $\dim \HHH=g$ and $\Bs
\HHH=\Phi^{-1}\Phi(o)$. In particular, $\HHH$ has no fixed
components. Moreover, the image of the map $\Phi_{\HHH}$ coincides
with the image of the projection of $\ov X\subset \PP^{g+1}$ from
the point $\Phi(o)$.

\begin{claim}
$\dim \Phi_{\HHH}(X)\ge 2$ and $\dim \Phi_{\HHH}(X)=2$ if and only
if $\ov X\subset \PP^{g+1}$ is a projective cone over the surface
$\Phi_{\HHH}(X)$
\end{claim}
\begin{proof}
Obvious.
\end{proof}

\begin{claim}
If $X$ is not hyperelliptic, then $\dim \Phi_{\HHH}(X)=3$.
\end{claim}
\begin{proof}
Assume that $\Phi$ is an isomorphism and $\dim \Phi_{\HHH}(X)=2$.
Then $\ov X$ is a projective cone over a K3 surface. But in this
case the vertex of the cone $X=\ov X$ is not a canonical
singularity (see, e.g., \cite[2.14]{Reid}), a contradiction.
\end{proof}

\begin{emptytheorem}
\label{pusto-proof}
Consider the case $\dim \Phi_{\HHH}(X)=3$. Let $f\colon (\hat
X,\hat\HHH)\to (X,\HHH)$ be a relative terminal model. Thus the
pair $(\hat X,\hat\HHH)$ has only terminal singularities and
\[
K_{\hat X}+\hat\HHH=f^*(K_X+\HHH)-E,
\]
where $E=\sum a_iE_i$ is an integral effective Weil divisor. Since
$(X,\HHH)$ is not canonical, $E\neq 0$ (see Theorem
\xref{th-can-sing}). Run $K_{\hat X}+\hat\HHH$-MMP: $(\hat
X,\hat\HHH)\dashrightarrow (W,\HHH_W)$. It is clear that $K_{\hat
X}+\hat\HHH\equiv -E$, so our MMP is the same as $-E$-MMP. At the
end we get a fiber type contraction $h\colon W\to Z$ which is
$E_W$-positive, where $E_W$ is the birational transform of $E$.
Indeed, $K_{W}+\HHH_W\equiv -E_W$ cannot be nef. Note that the
image $\Phi_{\HHH_W}(W)=\Phi_{\HHH}(X)$ is three-dimensional, so
$\HHH_W$ is not a pull-back of a linear system on $Z$. Thus
$\HHH_W$ is positive on the fibers of $h$. Since $(W,\HHH_W)$ is
terminal, the linear system $\HHH_W$ has at most nonsingular
isolated base points of multiplicity $1$ (see Proposition
\xref{prop-A-sing-term}). In particular, $\HHH_W$ is a linear
system of Cartier divisors. If $Z$ is not a point, then $h$ is
either a generically $\PP^1$, $\PP^2$, or $\PP^1\times
\PP^1$-bundle and $X$ is rational (see, e.g., \cite{A}). If $Z$ is
a point, then $W$ is a $\QQ$-Fano with $\rho=1$. Let $H\in \HHH_W$
be a general member. Then $H$ is a smooth surface and by the
Adjunction Formula $-K_H$ is ample, i.e., $H$ is a del Pezzo
surface. Assuming that $W$ is nonrational from Campana-Flenner's
Theorem \xref{th-CF} we get cases
\xref{gather-cases-1}-\xref{gather-cases-4}. So, $g=\dim \HHH\le
\dim |\HHH_W|\le 4$, where $|\HHH_W|$ is the complete linear
system generated by $\HHH_W$. Hence in case \xref{gather-cases-4}
we have $g=2$ and the variety $X$ is hyperelliptic. From classical
results on K3 surfaces we get that $X$ is either a quartic or a
complete intersection of a quadric and a cubic (cf. \cite{Isk}).
This proves Theorem \xref{rational} in the case $\dim
\Phi_{\HHH}(X)=3$.

If $\dim \Phi_{\HHH}(X)=2$, then we can use the same arguments.
The only difference is that $\HHH_W$ can be a pull-back of a
linear system on $Z$. This is possible when $\dim Z=2$, i.e., when
$h\colon W\to Z$ is a conic bundle. Then we get case
\xref{case-hyperell-conic}. Theorem \xref{rational} is proved.
\end{emptytheorem}

\begin{corollary}
In cases \xref{case-quartic} and \xref{case-quadric-cubic} the
variety $X$ is unirational.
\end{corollary}
\begin{proof}
Follows by the fact that all varieties in
\xref{gather-cases-1}-\xref{gather-cases-3} are unirational (see
\cite{CF}, \cite[Ch. 4]{Manin}, and also \cite[\S 10.1]{IP}).
\end{proof}

\begin{example}
\label{ex-quartic}
Consider the quartic threefold $X$ defined by
$x_0^2x_4^2+\phi(x_1,\dots,x_4)=0$ in $\PP^4$, where
$\phi(x_1,\dots,x_4)$ is a sufficiently general homogeneous
polynomial of degree $4$. By Bertini theorem the only singularity
of $X$ is the point $P=(1, 0, 0, 0, 0)$. It is easy to see that
$P$ is canonical and is not cDV (see \cite{Reid}). Therefore, $X$
satisfies conditions of Theorem \xref{rational}. Clearly, $X$ is
birationally equivalent to the weighted hypersurface
$x_0^2+\phi(x_1,\dots,x_4)=0$ in $\PP(2,1,1,1,1)$, i.e., a
nonsingular del Pezzo threefold of degree $2$. The last variety is
known to be nonrational (see \cite{Be}).

On the other hand, a quartic with a single triple point is
obviously rational.
\end{example}

We give a more complicated example of a rational quartic.

\begin{example}
Consider the quartic $X$ given by $x_0^2x_4^2+x_1^4+x_2^4+x_3^4=0$
in $\PP^4$. This $X$ has two non-cDV points. The map
\[
X \dashrightarrow\PP^3, \quad (x_0,x_1,x_2,x_3,x_4)
\dashrightarrow (x_1,x_2,x_3,x_4)
\]
is generically finite of degree $2$. This shows that $X$ is
birationally equivalent to del Pezzo threefold $X'$ of degree $2$
that is a double cover of $\PP^3$ branched along the quartic
$(x_1^4+x_2^4+x_3^4=0)\subset \PP^3$. Since $X'$ has a non-cDV
point, it is rational (see Lemma \xref{lemma-del-Pezzo} below).
\end{example}

\begin{example}
\label{ex-quardic-cubic}
Let $V\subset \PP^5$ be a cone over a cubic threefold $T$ having
at worst isolated double points, let $Q$ be a sufficiently general
quadric passing through the vertex $P$, and let $X=V\cap Q$. Then
$X$ satisfies conditions of Theorem \xref{rational}. The
projection from $P$ induces a birational isomorphism between $X$
and $T$. Thus $X$ is rational if and only if $T$ is smooth
\cite{CG}.
\end{example}

\section{Hyperelliptic case}
In this section we study hyperelliptic Fano threefold.
\subsection*{Notation}
Let $\EEE=\OOO_{\PP^1}(a_1)\oplus \OOO_{\PP^1}(a_2)\oplus
\OOO_{\PP^1}(a_3)$, where $a_1\ge a_2\ge a_3\ge 0$. Consider the
scroll $\WW=\WW(a_1,a_2,a_3)=\PP_{\PP^1}(\EEE)$. Let $M$ be the
tautological divisor and let $F$ be a fiber of the projection
$\pi\colon \WW\to \PP^1$. In this notation, $-K_{\WW}=3M+(2-\sum
a_i)F$. The linear system $|M|$ is base point free and defines a
morphism $\var\colon \WW\to \ov \WW\subset \PP^{\sum a_i+2}$. It
is easy to see that $\deg \ov \WW=\sum a_i$.

\begin{emptytheorem}
Let $X$ be a hyperelliptic Fano threefold of genus $g$ with at
worst canonical Gorenstein singularities and let $\Phi\colon X\to
\ov X\subset \PP^{g+1}$ be the anticanonical morphism. Let
$B\subset \ov X$ be the branch divisor. By the ramification
formula we have
\begin{equation}
\label{eq-Hurwitz}
K_X\sim \Phi^*\left(K_{\ov X}+\frac12B\right),
\end{equation}
where $\frac12B$ is a class of an integral Weil divisor.
Therefore, $-(K_{\ov X}+\frac12B)$ is linearly equivalent to a
hyperplane section $L$ of $\ov X\subset \PP^{g+1}$. Thus
\[
B\sim 2(-K_{\ov X}-L).
\]
\end{emptytheorem}
\begin{emptytheorem}
According to Proposition \xref{prop-hyperell} and Remark
\xref{rem-Enriques} we can distinguish the following
possibilities:
\begin{subth}
\label{cases-P3}
$g=2$, $\ov X=\PP^3$ and $B\sim
\OOO_{\PP^3}(6)$,
\end{subth}
\begin{subth}
$g=3$, $\ov X=Q\subset \PP^4$ is a (possibly singular) quadric,
$B\sim \OOO_{Q}(4)$,
\end{subth}
\begin{subth}
\label{cases-Veronese-cone}
$g=5$, $\ov X=V_4\subset \PP^6$ is a Veronese cone, $B\sim
\OOO_{V_4}(3)$,
\end{subth}

\begin{subth}
\label{cases-a,a,a}
$\ov X=\WW(a_1,a_2,a_3)$, where $a_1\ge a_2\ge a_3\ge 1$, $g=\sum
a_i +1\ge 4$, $B\sim 2(2M+(3-g)F)$,
\end{subth}

\begin{subth}
\label{cases-a,a,0}
$\ov X=\ov{\WW(a_1,a_2,0)}$, where $a_1\ge a_2\ge 1$, $g=a_1+a_2
+1\ge 4$, $B=\var(B')$, where $B'\sim 2(2M+(3-g)F)$,
\end{subth}

\begin{subth}
\label{cases-a,0,0}
$\ov X=\ov{\WW(g-1,0,0)}$, where $g\ge 4$, $B\sim 2(g+1)P$, where
$P$ is the class of a plane on $\ov{\WW}\subset \PP^{g+1}$.
\end{subth}
\end{emptytheorem}

\begin{remark}
Similar to \cite[Th. 2.2]{Isk} one can get a complete
classification of hyperelliptic Fano threefolds. Indeed, by
\eqref{eq-Hurwitz} the pair $(\ov X, \frac12 B)$ is klt and
discrepancies of $(\ov X, \frac12 B)$ are contained in
$\frac12\ZZ$. This fact gives us very strong restrictions on
$(a_1,a_2,a_3)$. For example, in case \xref{cases-a,0,0} the
singularity of $\ov X$ along the vertex of the cone is locally
isomorphic to $\CC^2/\ZZ_{g-1}(1,1)\times \CC$. Since
$\mt{discr}(\ov X)\ge\mt{discr}(\ov X, \frac12 B)\ge -1/2$, we
have $-1+2/(g-1)\ge -1/2$, so $g\le 5$. We get only two
possibilities: $\ov X$ is a weighted hypersurface $W_{10}\subset
\PP(1,1,3,3,5)$ or $W_{12}\subset \PP(1,1,4,4,6)$.
\end{remark}

The following proposition is very easy but the author could not
find a reference.

\begin{proposition}
Let $X$ be a Fano threefold with at worst canonical Gorenstein
singularities. If $X$ is hyperelliptic, then $X$ is unirational
except perhaps for cases \xref{cases-P3} and
\xref{cases-Veronese-cone}.
\end{proposition}
\begin{proof}
We consider only cases \xref{cases-a,a,a} and \xref{cases-a,a,0}.
Other cases are similar. Then $\var\colon \WW\to \ov X$ is either
an isomorphism or a small morphism (contracting a section). We
have the following diagram
\[
\begin{diagram}
X&\lTo^{\Psi}&Y&&\\
\dTo^{\Phi}&&\dTo^{\phi}&\rdTo^{\psi}&\\
\ov X&\lTo^{\var}&\WW&\rTo^{\pi}&\PP^1
\end{diagram}
\]
where $Y$ is the normalization of the dominant component of
$X\times_{\ov X} \WW$. The morphism $\Psi$ does not contract
divisors. Hence, $K_Y=\Psi^*K_X$ and $Y$ has only canonical
Gorenstein singularities.

Denote by $B'$ the birational transform of $B$ on $\WW$. Then
\[
B'\sim 2(2M+(3-g)F).
\]
Since $F\cap B'\neq \emptyset$ for any fiber $\pi^{-1}(\mt{pt})$,
the restriction $\phi\circ\phi^{-1}(F)\colon \phi^{-1}(F)\to F$ is
not \'etale. Hence all fibers $\phi^{-1}(F)=\psi^{-1}(\pi(F))$ are
connected. Now one can see that the general fiber $Y_\eta$ is a
del Pezzo surface with at worst Du Val singularities. Further,
\[
K_{Y_{\eta}}^2=K_Y^2\cdot Y_\eta= 2\left(K_{\WW}+\frac12
B'\right)^2\cdot F=2.
\]
Therefore, $\psi \colon Y\to \PP^1$ is a (possibly singular) del
Pezzo fibration of degree $2$. The  In this case $Y$ is
unirational because so is the general fiber over a non-closed
field $K(\PP^1)$ (see \cite[Ch. 4]{Manin}, \cite{CT}).
\end{proof}

\subsection*{Double Veronese cone}
Now we study Fano threefolds of type \xref{cases-Veronese-cone}.
We show that a nonrational Fano threefold from
\xref{case-hyperell-conic} or \xref{case-hyperell-other} cannot be
of this type. Note that in this case $X$ is isomorphic to a
weighted hypersurface $W_{6}\subset \PP(1,1,1,2,3)$ and
$-K_{X}\sim \OOO_{W_6}(2)$. Since $\OOO_{W_6}(1)$ is invertible,
$X$ is a del Pezzo threefold of degree $1$.
\begin{lemma}
\label{lemma-del-Pezzo}
Let $X$ be a del Pezzo threefold of degree $d$ \textup(with at
worst canonical Gorenstein singularities\textup). Assume either
$\rho(X)>1$ or singularities of $X$ are worse than terminal
factorial \textup(= isolated cDV\textup). Then $X$ is birationally
equivalent to $\PP^3$ or a del Pezzo threefold of degree $d'> d$
\textup(again with at worst canonical Gorenstein
singularities\textup). Moreover, if $X$ has a non-cDV point, then
$X$ is rational.
\end{lemma}
\begin{proof}
Consider $\QQ$-factorial terminal modification $f\colon Y \to X$.
Thus $Y$ has at worst factorial terminal singularities and
$-K_Y=2H$, where $H$ is a nef and big Cartier divisor. Let
$h\colon Y\to Z$ be a $K$-negative extremal contraction. By our
assumption $Z$ is not a point. If $Z$ is a curve or a surface,
then as in \xref{pusto-proof} we have that $h$ is either a
generically $\PP^1$, $\PP^2$, or $\PP^1\times \PP^1$-bundle. In
these cases $Y$ is rational.

Thus we may assume that $h$ is birational.  By the classification
of extremal rays on Gorenstein terminal threefolds \cite{Cut} we
have only one possibility:  $h$ contracts a divisor $E$ to a point
$P$. Let $H\in |H|$ be a general member and let $H_Z=h(H)$. We
claim that the point $P\in Z$ is smooth. Indeed, $H_Z$ is normal
outside of $P$. Since the point $P$ is terminal and
$\QQ$-factorial, $H_Z$ is normal also at $P$. On the other hand,
the contraction $h|_H\colon H\to H_Z$ is $K_H$-negative. In this
situation the point $\phi(H\cap E)\in H_Z$ must be smooth. This
implies that so is $P\in Z$. Again by \cite{Cut} $h$ is the
blow-up of $P$ and $K_{Y}=h^*K_Z+2E$. Thus $-K_Z=2H_Z$, where
$H_Z$ is a nef and big Cartier divisor, and $H^3_Z=H^3+1$. If
$\rho(Z)>1$, we can repeat procedure replacing $Y$ with $Z$. Since
$H^3\le 9$ the procedure terminates and we get a nonbirational
contraction. This proves the first part of the lemma.

To prove the second part we assume that $X$ has a non-cDV point,
say $P$. Then there is a two-dimensional component $S\subset
f^{-1}(P)$. For any curve $\Gamma\subset S$ we have $K_Y\cdot
\Gamma=0$. Therefore, $S\cap E=\emptyset$ and $h(S)$ again
satisfies the above property. This shows that after birational
contractions we cannot obtain a model with $\rho=1$. Hence $X$ is
rational.
\end{proof}

\begin{corollary}
Let $X$ be a hyperelliptic Fano variety of type
\xref{cases-Veronese-cone}. If $X$ has a non-cDV point, then it is
rational. If either $X$ has a nonterminal singularity or is not
$\QQ$-factorial, then it is unirational.
\end{corollary}

In conclusion we give an example of a nonrational hyperelliptic
Fano threefold such as in \xref{case-hyperell-conic}.
\begin{example}
Let $X$ be following weighted hypersurface in $\PP(1,1,1,1,3)$:
\[
x_0^2(x_1^4+x_2^4+x_3^4)+x_1^6+x_2^6+x_3^6=x_4^2.
\]
The only singularity of $X$ is at the point $P=(1, 0, 0, 0, 0)$.
In the affine chart $(x_0\neq 0)$ this singularity is given by
\begin{equation}
\label{eq-ex-hyp-conic}
x_1^4+x_2^4+x_3^4+x_1^6+x_2^6+x_3^6=x_4^2.
\end{equation}
It is easy to see that $P\in X$ is a hypersurface canonical
non-cDV singularity. By Adjunction $-K_X\sim \OOO_{\PP}(1)$ and
$-K_X^3=2$. Therefore, $X$ is a hyperelliptic Fano threefold of
genus $2$. The projection
\[
X\subset \PP(1,1,1,1,3)\dashrightarrow \PP^2,\qquad (x_0, x_1,
x_2, x_3, x_4)\dashrightarrow (x_1, x_2, x_3)
\]
gives us a structure of fibration into rational curves. Now let
$f\colon \hat X\to X$ be the weighted blowup of $P$ with weights
$(1,1,1,2)$ in the affine chart $(x_0\neq 0)$ (see
\eqref{eq-ex-hyp-conic}). Then $\hat X$ is smooth and (in notation
of \xref{pusto-proof}) $\hat \HHH$ is base point free. We get a
conic bundle $h\colon \hat X=W\to \PP^2=Z$. The discriminant curve
is given by the equation $(x_1^4+x_2^4+x_3^4)
(x_1^6+x_2^6+x_3^6)=0$ on $\PP^2$. Therefore, $X$ is nonrational
(see \cite{Sho1}).
\end{example}

\end{document}